\RequirePackage{amsthm}
\documentclass[pdflatex,sn-aps]{sn-jnl}[11pt]
\usepackage{amsmath, amsthm, amssymb}
\usepackage{lipsum,mathtools,graphicx}
\usepackage[dvipsnames]{xcolor}
\usepackage{xspace,tikz-cd}
\usepackage[title]{appendix}%
\usepackage[smalltableaux]{ytableau}
\usepackage{titlesec}
\titleformat{\section}
  {\normalfont\fontsize{12}{15}\bfseries}{\thesection}{1em}{}

\theoremstyle{thmstyleone}
\newtheorem{thm}{Theorem}[section]

\newtheorem{proposition}[thm]{Proposition}
\theoremstyle{thmstyletwo}
\newtheorem{conjecture}{\textbf{Conjecture}}
\theoremstyle{thmstylethree}
\newtheorem{example}[thm]{Example}
\theoremstyle{definition}
\newtheorem{definition}[thm]{Definition}
\newtheorem{remark}[thm]{Remark}
\newcommand{\qw}{$q$-Whittaker\xspace}
\newcommand{\mhl}{modified\,Hall-Littlewood\xspace}

\DeclareMathOperator{\ch}{ch}
\DeclareMathOperator{\cch}{cch}
\DeclareMathOperator{\inv}{inv}
\DeclareMathOperator{\quinv}{quinv}
\DeclareMathOperator{\maj}{maj}
\DeclareMathOperator{\dg}{dg}

\DeclareMathOperator{\Des}{Des}

\DeclareMathOperator{\leg}{leg}

\DeclareMathOperator{\imf}{IMF}
\DeclareMathOperator{\izf}{IZF}
\DeclareMathOperator{\qmf}{QMF}
\DeclareMathOperator{\qzf}{QZF}
\usepackage{hyperref}
\hypersetup{
    colorlinks=true,
    linkcolor=blue,
    filecolor=magenta,      
    urlcolor=cyan,
    }
\usepackage[dvipsnames]{xcolor}
\colorlet{Mycolor1}{green!80!blue}
\colorlet{Mycolor6}{lime!80!Mahogany}
\colorlet{Mycolor3}{SkyBlue!90!yellow}
\colorlet{Mycolor7}{Goldenrod!30!Orchid}
\colorlet{Mycolor5}{orange!80!RubineRed}
\colorlet{Mycolor4}{blue!40!pink}
\colorlet{Mycolor2}{red!80!RubineRed}

\begin{document}
\title[Bijections between different combinatorial models ]{Bijections between different combinatorial models for $q$-Whittaker and modified Hall-Littlewood polynomials}

\author{\fnm{T V} \sur{RATHEESH}}\email{ratheeshtv@imsc.res.in}

\affil{\orgdiv{ The Institute of Mathematical Sciences}, \orgname{A CI of Homi Bhabha National
Institute}, \orgaddress{\street{ Chennai}, \postcode{600113}, \country{India}}}
\abstract
{ We consider the monomial expansion of the $q$-Whittaker polynomials and the \mhl polynomials arising from specialization of
the modified Macdonald polynomial. The two combinatorial formulas for the
latter, due to Haglund, Haiman, and Loehr and Ayyer,
Mandelshtam and Martin, give rise to two different parameterizing sets in
each case. We produce bijections between the parameterizing sets, which
preserve the content and major index statistics.
We identify the major index with the charge or cocharge of appropriate
words, and use descriptions of the latter due to  Lascoux--Sch\"{u}tzenberger and Killpatrick to show that our bijections have the desired properties.}

\keywords{$q$-Whittaker polynomials, \mhl polynomials, charge and cocharge of a word, inv, quinv, and maj of a filling.}
\maketitle
\section{Introduction}
Let $x=(x_1, x_2,\ldots,)$ and $q, t$ denote a collection of algebraically independent indeterminates. The modified Macdonald polynomials $\widetilde{H}_\lambda(x;q,t)$ defined  by Garsia and Haiman \cite{garsia1993graded}, are  an important class of symmetric functions in the variables $x_i$, indexed by partitions $\lambda$. From the the well-known identity $\widetilde{H}_\lambda(x;q,t) =
\widetilde{H}_{\lambda'}(x;t,q)$ \cite[Prop 2.6]{haiman1999macdonald}, and by interchanging the role of $t$ and $q$ in \cite[(3.1)]{bergeron2020survey}, we can expand $\widetilde{H}_\lambda(x;q,t)$ in powers of $q$:
\begin{equation} 
\label{eq1.1}
\widetilde{H}_\lambda(x; q,t) =  Q_{\lambda^{\prime}}(x; t) q^{n(\lambda^{\prime})}  + \ldots +\widetilde{H}_{\lambda}(x; t),
\end{equation}
where $\lambda^{\prime}$ denotes the conjugate of $\lambda$, and $n(\lambda) =\displaystyle\sum_{k=0}^{\lambda_{1}}\binom{{\lambda_k}'}{2}$.
The leading coefficient  $Q_{\lambda^{\prime}}(x; t)$ is the so-called \qw polynomial \cite{garsia1993graded,bergeron2020survey} and the constant term  $\widetilde{H}_{\lambda}(x; t)$ is the  \mhl polynomial \cite[Section 2]{mellit2020poincare}.

In \cite{HHL01}, Haglund, Haiman, and Loehr provided a combinatorial interpretation of  $\widetilde{H}_\lambda(x;q,t)$. They considered fillings of the Young diagram of $\lambda$. To each filling, they associated an $x$-weight (a monomial in the $x_i$) along with $q$- and $t$-weights in terms of statistics called $\inv$ and $\maj$. In these terms, the $\widetilde{H}_\lambda(x;q,t)$ is just the generating function of all fillings counted with these weights. Later in \cite{ALCO_2023__6_1_243_0} Ayyer, Mandelshtam and, Martin found another statistic called $\quinv$ and showed that the formula for $\widetilde{H}_\lambda(x;q,t)$ still holds if $\inv$ is replaced with $\quinv$. This led them to conjecture that there should exist a ``natural'' bijection from the collection of all fillings to itself, which preserves both $x$-weight and $\maj$ and takes $\inv$ to $\quinv$.

In other words, for each $0 \leq p \leq n(\lambda^{\prime})$, one seeks to describe a bijection from the collection of all fillings having $\inv = p$  to the collection of all fillings having $\quinv=p$, which preserves both $\maj$ and $x$ weight. The main results of this paper are constructions of such bijections for the extreme values of $p$, i.e., for $p = n(\lambda^{\prime})$ and $p = 0 $.

To this end, we first show that on the set of fillings which maximise (resp. minimise, i.e., make zero)  $\inv$ or $\quinv$, the $\maj$ statistic coincides with the Lascoux--Sch\"{u}tzenberger charge (resp. cocharge) statistic on a suitable word. We use two descriptions of charge - the classical one due to Lascoux--Sch\"{u}tzenberger and a variant due to Killpatrick - to show that the natural bijection we construct between the $\inv$ and $\quinv$ pictures preserves $\maj$.

\section*{Acknowledgement}
The author would like to thank Aritra Bhattacharya and S Viswanath for helpful discussions and suggestions. The author acknowledges partial support under a DAE Apex Grant to the Institute of
Mathematical Sciences.

\section{Preliminaries}
Let $ \lambda = (\lambda_{1} \geq \lambda_{2} \geq \cdots)$ be a partition, and let $ \dg(\lambda) $ be its Young diagram, which we draw in English convention, as a left-up justified array of boxes with $\lambda_i$ boxes in the $i^{th}$ row from the top. Denote the box in the $i^{th}$ row and $j^{th}$ column  as the cell $u = (i ,j)$.

A {\em filling} of $ \dg(\lambda) $ is a function $ \sigma : \dg(\lambda) \to \mathbb{N}$, i.e., it assigns a positive integer to each cell in $\dg(\lambda)$. Let $\mathcal{F}(\lambda)$  be the set of fillings of $ \dg(\lambda) $, and let $|\lambda|= \displaystyle\sum_{i} \lambda_i$.

A {\em descent}  of a filling $\sigma$ is a pair of entries $\sigma(u) > \sigma(v),$ where $u, v \in \dg(\lambda)$ and the cell $u$ is immediately below  $v$, that is if  $ v = (i, j) $ then $ u = (i+1, j)$.
 Define
\begin{align}
\label{eq 1.0.5}
 \Des(\sigma) = \; \{u \in \dg(\lambda) : \sigma(u) > \sigma(v) \text{ is a descent}\}
\end{align}
and
\begin{align}
\maj(\sigma) = \sum_{u \in \Des(\sigma)} (\leg(u)+1),
\end{align}
where $\leg(u)$ is the number of cells strictly below $u$ in its column.

\begin{example} \label{Eg.1}
    For $n=8$ and $\lambda =(6, 4, 2)$, let 
\begin{align} \notag
 \sigma = \begin{ytableau}
       1 & 3 & 5 & 5 & 6 & 8\\
       2 & 4 & 2 & 7\\
       3 & 1 \\
       \end{ytableau}
\end{align}
    then  $\Des(\sigma) = \{(2, 1), (3, 1), (2, 2), (4, 2)\} 
      \text{ and } \maj(\sigma) = 6$.
    \end{example}
    
Given $\sigma$ in $\mathcal{F}(\lambda)$,
\begin{enumerate} 
  \item  The {\em inversion type augmented filling} $\hat{\sigma}$ is obtained by adding a cell with entry $0$ above the first cell in each column of $\sigma$.
  \item The {\em $q$-inversion type augmented filling} $\widetilde{\sigma} $ is obtained by adding a cell with entry $\infty$ below the last cell  in each column of $\sigma$.
\end{enumerate}
We denote the shapes of the resulting diagrams by $\dg(\hat{\sigma})$ and $\dg(\widetilde{\sigma})$ respectively. 
An {\em  inversion triple} in $\sigma $ is a triple of cells $(u, v, w)$ in $\dg(\hat{\sigma})$ such that:
\begin{enumerate}
    \item $u, w \in \dg(\lambda)$, and $w$ is in the same row as $u$, to its right,
    \item $v$ is immediately above $u$, in its column, and 
    \item the triple of values $(\hat{\sigma}(u), \hat{\sigma}(v), \hat{\sigma}(w))$, belongs to the set $\mathcal{I}$, where 
      \begin{equation}\label{eq:inv-trip}
        \mathcal{I}= \{ (x, y, z ) : y < z < x,\; x < y < z,\; z < x < y, \text{ or } x =  y \neq z \}.
      \end{equation}
\end{enumerate} 
Let $\inv(\sigma)$ denote the total number of {\em  inversion triples} in $\sigma$.
Likewise, a {\em  $q$-inversion triple} in $\sigma$ is a triple of cells $(u^{\prime}, v^{\prime}, w^{\prime})$ in $\dg(\widetilde{\sigma})$ such that: 
\begin{enumerate}
    \item $v^{\prime}, w^{\prime} \in \dg(\lambda)$ and $w^{\prime}$ is in the same row as $v^{\prime}$, to its right,
    \item $u^{\prime}$ is immediately below $v^{\prime}$, in its column, and
    \item The triple of values $(\widetilde{\sigma}(u^{\prime}), \widetilde{\sigma}(v^{\prime}), \widetilde{\sigma}(w^{\prime}))$, belongs to the set $\mathcal{Q}$, where 
      \begin{equation}\label{eq:quinv-trip}
            \mathcal{Q} = \{ (x, y, z) : y < z < x,\; x < y < z, \: z < x < y  \text{ or } x =  y \neq z \}.
          \end{equation}
\end{enumerate} 
We let $\quinv(\sigma)$ denote the total number of {\em  $q$-inversion triples} in $\sigma$.

For a more detailed discussion of the $\inv$ and $\quinv$ statistics, we refer the reader to \cite{HHL01} and \cite{ALCO_2023__6_1_243_0}.
\begin{example}
For $\sigma$ from Example \ref{Eg.1}, 
\begin{align} \notag
 \hat{\sigma} = \begin{ytableau}
       0&0&0&0&0&0\\
       1 &3 & 5 & 2 & 6 & 8\\
       2 &4 & 1 & 7\\
       3 & 1 \\
\end{ytableau} 
\end{align}
and
\begin{align} \notag
 \widetilde{\sigma} = \begin{ytableau}
       1 &3 & 5 & 2 & 6 & 8\\
       2 &4 & 1 & 7 &\infty&\infty\\
       3 & 1 &\infty&\infty\\
       \infty&\infty\\
\end{ytableau}\\ \notag
    \end{align}
     with  $\inv(\sigma) = 3 $, and $\quinv(\sigma) = 6$. 
\end{example}

For a filling $\sigma$, its {\em content} is the sequence $\mu = (\mu_{1}, \mu_{2}\ldots , )$ of non negative integers, where $\mu_{i} = |\{u \in \dg(\lambda): \sigma(u) = i \}|$ is the number of $i$'s in $\sigma$. Let    \begin{equation}
x^{\sigma}= \prod_{i \geq 1 } x^{\mu_{i}}\, . \nonumber
\end{equation}
The following are the main theorems of \cite{HHL01} and \cite{ALCO_2023__6_1_243_0}.
\begin{thm}
\label{Theorem 1}
\cite[Theorem 2.2] {HHL01} Let $\lambda$ be a partition. The modified Macdonald polynomial can be written as
\begin{equation}
\label{eq1.0.9}
\widetilde{H}_\lambda(x;q,t) = \sum_{\sigma \in \mathcal{F}(\lambda)} q^{\inv(\sigma)} t^{\maj(\sigma)} x^{\sigma}.
\end{equation}
\end{thm} 

\begin{thm} \label{Theorem 2}
\cite[Theorem 2.6]{ALCO_2023__6_1_243_0}
 Let $\lambda$ be a partition. The modified Macdonald polynomial can be written
as
\begin{equation}
\label{eq1.0.13}
\widetilde{H}_\lambda(x; q,t) = \sum_{\sigma \in \mathcal{F}(\lambda)} q^{\quinv(\sigma)}t^{\maj(\sigma)}x^{\sigma}.
\end{equation}
\end{thm}

From the above theorems one can naturally ask whether there exists a bijection on  $\mathcal{F}(\lambda)$, that preserves both $\maj$ and content, and takes $\inv$ to $\quinv$. In  \cite{ALCO_2023__6_1_243_0}, the authors conjectured an even stronger statement, which we explain below.

\begin{definition}
  Let $\sigma, \tau$ in $\mathcal{F}(\lambda)$, we say $\sigma$ and $\tau$  are {\em row-equivalent} if every row of $\sigma$ is a permutation of the corresponding row of $\tau$. We write $\sigma \sim \tau$  when $\sigma$ and $\tau$ are row-equivalent.
\end{definition}
It will also be convenient to define $R_i(\sigma)$ to be the multiset of entries in the $i^{th}$ row of $\sigma$. In terms of this notation,  $\sigma \sim \tau$ iff $R_i(\sigma) = R_i(\tau)$ for all $i$.

\begin{conjecture} (\cite{ALCO_2023__6_1_243_0}):
 Given a partition $\lambda$, there exists a bijection $\delta :\mathcal{F}(\lambda) \to \mathcal{F}(\lambda)$
such that for all $\sigma$,
\begin{align} \notag
    \sigma &\sim  \delta(\sigma),\\ \notag
    \maj(\delta(\sigma)) &=  \maj(\sigma),\\ \notag
    \quinv(\sigma) &= \inv(\delta(\sigma)).
\end{align}
\end{conjecture}
In the following sections, we prove the conjecture for the extreme values of $\inv$ and $\quinv$\footnote{The conjecture in \cite{ALCO_2023__6_1_243_0} is stated for fillings with an upper bound stipulated for the entries, but the version
stated here can be seen to be equivalent to that.}.

\subsection{Two descriptions of the charge and cocharge statistics}
In this subsection we will recall the classical definition of the statistics {\em charge}  and {\em cocharge} on words, introduced by Lascoux and Sch\"{u}tzenberger \cite{lascoux1978conjecture} and a new description to calculate them due to  Killpatrick \cite{killpatrick2000combinatorial}. It will emerge that the classical definition is more attuned to the  $\inv$ statistic and the one due to Killpatrick corresponds to the $\quinv$ statistic.

\subsubsection{Classical description of charge}\label{sec:ls-charge}
Let $w$ be a (finite) word with letters in the alphabet $\mathbb{N}$. The content of $w$ is the sequence $(\mu_i)_{i \geq 1}$ where $\mu_{i}$ is the number of occurrences of $i$ in $w$. We say that $w$ is of {\em partition content} if $\mu_{i} \geq \mu_{j}$  for all $i < j $. The length of a word is defined as the total number of letters in it.
A word  with content $(1,1,\ldots , 1,0,0,\ldots)$ is called a  {\em standard word}. 
\begin{definition}
   Let $w$ be a standard word $w$ of length $k$. The {\em charge} of $w$, denoted $\ch(w)$, is defined as $\displaystyle\sum_{i} (k-i)$, where the sum ranges over all $1 \leq i \leq k$ for which  $i + 1$ occurs to the right of $i$ in $w$. 
\end{definition}

\noindent
For instance, for the standard word $w = 1465327$ of length $7$, $ \ch(w) =  6 + 3 + 1 = 10 $.
For a word $w$ with an arbitrary {partition content} $\mu =(\mu_{1}, \mu_{2}, \ldots)$, the {charge} $\ch(w)$ can be calculated as the sum of the charges of some particular standard subwords of $w$. The algorithm to choose such subwords proceeds as follows:

For the first standard subword, we begin from the right end of $w$, scan to the left until we come across a $1$. Then starting from that point, we scan left again until we encounter a $2$ and so on. Suppose for an  entry $j$ if there is no $j + 1$  to the left of it, then we wrap-around and continue to scan from the right end (i.e., picking the rightmost entry $j + 1$ in $w$). We continue this process until we reach the largest letter in $w$. The subset of the letters of $w$ thus chosen is now viewed as a subword of $w$, denoted $w_1$.

Once $w_{1}$ has been picked, we remove its letters from $w$ to obtain a new word $w'$. We repeat the above algorithm on $w'$ to obtain another standard subword $w_{2}$ of $w$. Continuing this process, we will get  $\mu_{1}$ standard subwords $w_{1}, w_{2}, \ldots, w_{\mu_{1}}$ of $w$. Define the {\em charge} of $w$ to be:
\[\ch(w) = \sum_{i=1}^{\mu_{1}}\ch(w_{i})\, .\]

\begin{example}
\label{eg_ls_1}
    Let $w = 121123132213$  be a word with content $(5,4,3)$. Then, the standard subwords of $w$ given by the above algorithm are $w_{1}= 321$, $w_{2}= 213$, $w_{3}= 213,\, w_{4}= 12,\text{ and } w_{5}= 1$  with the respective charges $\ch(w_{1}) = 0,\,\ch(w_{2}) = 1,\,\ch(w_{3}) = 1,\,\ch(w_{4}) = 1$ and  $\ch(w_{5}) = 0$. Therefore the charge of $w$ is $\ch(w) = 3. $
\end{example}

\subsubsection{New description of charge}\label{sec:killp} We now describe a variant of this procedure due to Killpatrick, which will play a key role in our bijections.
In this description, instead of the $w_{i}$'s we choose different standard subwords $v_{i}$'s as follows: 
\begin{enumerate}
    \item For the first standard subword, we begin from the left end of $w$, scan to the right until we come across the largest letter of $w$ (say $m$). Then, starting from that point, scan right again until we encounter an $m-1$, and continue in this order.
    \item If there is no $j$  to the right of the previous entry $j+1$, then wrap-around and continue scanning from the left end (i.e., pick the leftmost entry $j$). We continue this process until we  reach $1$. The letters in $w$ chosen as described are viewed as forming a subword $v_{1}$ of $w$. 
\item 
  Once $v_{1}$ is chosen, delete its entries from $w$ to obtain a new word $v'$. Repeat the above algorithm on the word $v'$  (starting with the largest letter that occurs in $v'$) to pick the standard subword  $v_{2}$ etc. Continue this process till we get  $\mu_{1}$ standard subwords $v_{1},\,v_{2},\, \ldots ,\, v_{\mu_{1}}$ of $w$.
\end{enumerate}

Define 
\[\widehat{\ch}(w) = \sum_{i=1}^{\mu_{1}}\ch(v_{i}) \,,\]
where $\ch(v_{i})$ is the {charge} of the standard word $v_{i}$ as defined in section \ref{sec:ls-charge}.

\begin{thm}\cite[Theorem 8] {killpatrick2000combinatorial}\label{Theorem_kill}
The statistic $\widehat{\ch}$ coincides with charge, that is, $\widehat{\ch} = \ch $.
\end{thm}

\begin{example}
    For $w$ from Example \ref{eg_ls_1}, the standard subwords of $w$ given by the new description are $v_{1}= 321$, $v_{2}= 132$, $v_{3}= 213$, $v_{4}= 21$, and $v_{5}= 1$ with respective charges  $\ch(v_{1}) = 0,\,\ch(v_{2}) = 2,\,\ch(v_{3}) = 1,\,\ch(v_{4}) = 0$  and $\ch(v_{5}) = 0$.
\end{example}

\subsubsection{Cocharge of a word}
\begin{definition}
For a standard word $w$ of length $k$, the {\em cocharge} of $w$, denoted $\cch(w)$, is defined as $\sum_{i}(k-i)$, where the sum ranges over all $i$, for which $i + 1$ occurs to the left of $i$ in $w$.
The {cocharge} of an arbitrary word $w$ with {partition content} $\mu =(\mu_{1}, \mu_{2}, \ldots)$, is defined as
\begin{align}\notag
     \cch(w) = &  \sum_{i=1}^{\mu_{1}}\cch(w_{i})\,,
\end{align}
where $w_{i}$'s are the standard subwords of $w$ chosen as in Section \ref{sec:ls-charge}.
\end{definition}

In \cite[Section 2.6]{butler1994subgroup}, an alternative definition for cocharge is provided in terms of charge. For a word $w$ with {partition content} $(\mu_i)_{i \geq 1}$, its {cocharge}  
\begin{align}
   \cch(w) &= \sum_{k= 1}^{\mu_{1}}\binom{\mu_{k}'}{2} - \ch(w) \, ,
\end{align}
where $\mu'$ is the conjugate of $\mu$. Thus, from Theorem \ref{Theorem_kill}, we have
\begin{align}\notag
     \cch(w) = &  \sum_{k=1}^{\mu_{1}}\big[\binom{\mu_{k}'}{2}-\ch(v_{k})\big]\, ,
\end{align}
where $v_{k}$'s are chosen as in Section \ref{sec:killp}. Since the length of $v_{k}$ is $\mu_{k}'$, we have $\cch(v_{k}) = \displaystyle\binom{\mu_{k}'}{2} -\ch(v_{k})$. This allows us to define the {cocharge} of the word $w$ as the sum of {cocharges} of the subwords, chosen as in Section \ref{sec:killp}, that is
\begin{align}
\label{eq.cch}
     \cch(w) = &  \sum_{k=1}^{\mu_{1}}\cch(v_{k}) \, .
\end{align}
Throughout this paper we use \eqref{eq.cch} as the definition of {cocharge}. 

For a more detailed discussion of the {charge} and {cocharge} statistics, we refer the interested reader  to  \cite[Chapter 2]{butler1994subgroup}.

\subsection{Charge and cocharge words}\label{sec:order}
Define the following total order on $\dg(\lambda)$:
  for two cells $u=(i,j)$ and $v= (i^{\prime},j^{\prime})$, we say  $u < v $, if:
  \begin{enumerate}
      \item  $i < i^{\prime}$,   or
      \item  $i= i^{\prime}$ and $j > j^{\prime}$.
  \end{enumerate}
  Given $\sigma$ in $\mathcal{F}(\lambda)$ or in $\mathcal{F}(\lambda,\,n)$, arrange the cells of $\dg(\lambda)$ in sequence: $u_{1} = (i_{1}, j_{1}), u_{2} = (i_{2}, j_{2}), \ldots , u_{k}=(i_{k}, j_{k})$,  where $\sigma(u_{1}) \geq \sigma(u_{2}) \geq  \ldots \geq \sigma(u_{k}) $ and $|\lambda| = k$. In cases where there is a constant segment $\sigma(u_{p}) = \sigma(u_{p+1}) =  \ldots = \sigma(u_{q}) $, the cells $u_{p},  u_{p+1} , \ldots , u_{q} $ are arranged in increasing order with respect to the total order $<$ defined above.

\begin{definition}(\cite[Section 7]{HHL01}) \label{def 2.2.1}
    The {\em cocharge word } of $\sigma$, denoted $cw(\sigma)$, is the word $ i_1 i _2 \ldots  i_k $ of content $\lambda$,  where $i_p$ is the row index of the cell $ u_p$ in the above order. The {\em charge word } $w(\sigma)$ is the reverse of  $cw(\sigma)$, that is $w(\sigma) = i_k i _{k-1} \ldots  i_1 $. 
\end{definition}

Next, instead of the total order $<$ on $\dg(\lambda)$, we consider a slightly different total order $<'$  defined as follows: for two cells $u=(i, j)$ and $v= (i^{\prime}, j^{\prime})$ we say  $u <' v $, if
  \begin{enumerate}
      \item  $i < i^{\prime}$,   or 
      \item  $i= i^{\prime}$ and $j < j^{\prime}\, .$
  \end{enumerate}
  For a filling $\sigma$, as above one can arrange the cells of $\dg(\lambda)$ in sequence with respect to the total order $<'$, and let $cw'(\sigma)$ (resp. $w'(\sigma)$) be the {cocharge word} (resp. {charge word}) of $\sigma$ with respect to $<'$, defined in the same way as in Definition \ref{def 2.2.1}. Then, it is easy to see that $cw'(\sigma) = cw(\sigma) $  and $w'(\sigma)= w(\sigma).$ 

 \begin{example}
    For $\sigma$ from Example \ref{Eg.1}, the ordering of cells of $\dg(\lambda)$ with respect to $\sigma$ and $<$ is:
$u_{1} = (1,6)$, $u_{2} = (2, 4)$, $u_{3} = (1, 5)$, $u_{4} = (1, 4)$, $u_{5} = (1, 3)$, $u_{6}= (2, 2)$, $u_{7} = (1, 2)$, $u_{8}= (3, 1 )$, $u_{9} = (2, 3)$, $u_{10} = (2, 1)$, $u_{11}= (1, 1)$ and $u_{12}= (3, 2)$.

On the other hand, the ordering with respect to $<'$ is: 
$v_{1} = (1,6)$, $v_{2} = (2, 4)$, $v_{3} = (1, 5)$, $v_{4} = (1, 3)$, $v_{5} = (1, 4)$, $v_{6}= (2, 2)$, $v_{7} = (1, 2)$, $v_{8}= (3, 1 )$, $v_{9} = (2, 1)$, $v_{10} = (2, 3)$, $v_{11}= (1, 1)$ and $v_{12}= (3, 2)$.

The corresponding words are $cw( \sigma) = cw'( \sigma) = 121112132213$ and $w( \sigma) = w'( \sigma)= 312231211121$.
\end{example}

\section{Bijection between the models for \texorpdfstring{$q-$}-Whittaker  polynomials} From Theorems \ref{Theorem 1}, \ref{Theorem 2}, and \eqref{eq1.1}, we may express the \qw polynomial as:
\begin{align} \notag
Q_{\lambda^{\prime}}(x;t) = &\sum_{\substack{\sigma \in \mathcal{F}(\lambda) \\ \inv(\sigma)= n(\lambda^{\prime})}}t^{\maj(\sigma)} x^{\sigma} \, ,\\ \notag
 = &\sum_{\substack{\sigma \in \mathcal{F}(\lambda) \\ \quinv(\sigma)=n(\lambda^{\prime})}}t^{\maj(\sigma)} x^{\sigma}\, . 
\end{align}
The following two propositions characterise the class of $\inv$ and $\quinv$ maximum fillings of $\dg(\lambda)$.

\subsection{inv maximisers}
\begin{proposition} \label{lemma 1} Let $\lambda=(\lambda_1 \geq \lambda_2 \geq \cdots)$  be a partition. For each $i \geq 1$, let $ S_{i}$ be any set of positive integers with $|S_i|= \lambda_{i}$. Then there exists a unique filling $\sigma \in \mathcal{F}(\lambda)$ with $\inv(\sigma) = n(\lambda^{\prime})$ such that the (multi)set $R_i(\sigma)$ of entries in the $i^{th}$ row of $\sigma$ is $S_{i}$ for all $i$. Furthermore, $\sigma$ satisfies $\maj(\sigma) = \ch(w({\sigma}))$.
\end{proposition}
\noindent
{\em Proof:} 
Let $\sigma$ be a filling with $\inv(\sigma) = n(\lambda ^{\prime})$. Then from the definition of $\inv$ in \eqref{eq:inv-trip}, and the maximality of $\inv(\sigma)$ one can observe the following facts about $\sigma$:
\begin{enumerate}
\item For each $i$, the entries in the $i^{th}$ row of $\sigma$ must be distinct
    \item $\sigma$ is strictly decreasing in row 1, and 
    \item for every cell $u$ not in row 1, if $v$ is the cell directly above $u$, and $S$ is the set
consisting of $u$ and the cells to its right in the same row, then 
\begin{enumerate}
    \item  If $ \sigma(u) \leq \sigma(v)$, then
  $\sigma(u) = \max \{ x\in \sigma(S) : x \leq \sigma(v)\}$, 
\item If $\sigma(u) > \sigma(v)$ then
$ x > \sigma(v)$ for all $x \in \sigma(S)$, and $\sigma(u) = \max \sigma(S)$.
\end{enumerate}
\end{enumerate}

From the above facts we can conclude that for a given collection $\{S_{i} : S_{i}\subset \mathbb{N}, |S_{i}| = \lambda_{i} \; \forall \; i \}$, the filling $\sigma$ with the required properties is unique if it exists. To establish existence, we construct $\sigma$ as follows:
\begin{enumerate}
    \item Take the entries of $\sigma$ in row 1 (from left to right) to be the elements of $S_1$ in strictly decreasing order.
    \item Inductively, once we have constructed rows $1$ through $i - 1$, the entries $\sigma(u)$ in row $i$ are determined one-by-one, moving  from left to right, as follows: (a) let $v$ be the cell directly above $u$ and define $\sigma(u)$ to be the largest unused element $x$ in $S_{i}$ such that $x \leq \sigma(v)$, provided such an $x$ exists (b) if no such $x$ exists, define $\sigma(u)$ to be the largest unused element in $ S_{i}$.
    \end{enumerate}
It is straightforward to verify that $\sigma$ satisfies the conditions in the first part of the proposition.

To prove $\maj(\sigma)=\ch(w({\sigma}))$, let $w_{1}$ be the first standard subword of $w({\sigma})$ in the classical definition of charge (Section \ref{sec:ls-charge}), and for all $i$, let $t_{i, 1}$ denote the letter $i$ in $w_{1}$. From the definition, $t_{1, 1}$ is the rightmost $1$ in $w({\sigma})$, and from the construction of $w({\sigma})$ (Section \ref{sec:order}), $t_{1,{1}}$ corresponds to the largest entry in the first row of $\sigma$, that is $t_{1,{1}}$ corresponds to $\sigma(1, 1)$. For $i \geq 1$, $t_{i+1,\,{1}}$ will be the rightmost $(i+1)$ to the left of $t_{i,{1}}$ in $w({\sigma})$ if $\sigma(i, 1) \geq \sigma(i+1, 1)$, and $t_{(i+1),{1}}$ will be the rightmost $i+1$ in $w({\sigma})$ otherwise. Observe that $t_{(i+1),\,{1}}$ lies to the right of $t_{i,{1}}$ in $w_{1}$ if and only if the pair  $\sigma(i+1, 1) > \sigma(i, 1)$ is a descent.
Therefore 
\[\ch(w_{1}) = \displaystyle\sum_{u \in \Des(C_{1}(\sigma))}(\leg(u)+1) \, ,\] where $\Des(C_{1}(\sigma))$ denotes the collection of cells $u$ where the pair $\sigma(u) >\sigma(v)$ is a descent in the first column of $\sigma$.

If we remove the first column of $\sigma$, the remaining  will be a filling $\sigma^{\prime}$ of $\dg(\mu)$, where $\mu$ is the shape of $\sigma'$. It is also easy to see that  $\sigma^{\prime}$ is an $\inv$-maximal filling of  $\dg(\mu)$, and the word $w({\sigma^{\prime}})$ corresponding to $\sigma^{\prime}$ is nothing but the subword of $w(\sigma)$ obtained by removing $w_{1}$.

Therefore, we can conclude that
\[\ch(w({\sigma})) = \ch(w_{1}) + \ch(w({\sigma^{\prime}}))\, ,\]
and 
\[ \maj(\sigma) =  \maj(\sigma') + \displaystyle\sum_{u \in \Des(C_{1}(\sigma))}(\leg(u)+1) \, .\]
It follows by induction that
$\ch(w({\sigma})) =  \maj(\sigma)$. \qed

\begin{remark}
The unique $\inv$-maximal filling corresponding to a given collection of sets is exactly the transpose of the unique {\em coinversion} zero non-attacking filling 
corresponding to that same collection of sets (see \cite[Proposition 17]{alexandersson2020cyclic}).
 \end{remark}
\begin{example}
\label{eg_invm}
    Let $\lambda = (7, 5, 4, 2) $, and let $S_{1}= \{1, 2, 3, 4, 5, 6, 7\}$, $S_{2}= \{2, 5, 7, 9, 10\}$, $S_{3}= \{6, 8, 9, 10\}$, and $S_{4}= \{7, 8\}$, then the corresponding $\inv$-maximal filling is $$\sigma = \begin{ytableau}
       *(Mycolor1)7 &*(Mycolor2) 6 & *(Mycolor3)5 & *(Mycolor4)4 &*(Mycolor5) 3 &*(Mycolor6)2 & *(Mycolor7)1 \\
     *(Mycolor1)  7 & *(Mycolor2)5 & *(Mycolor3)2 &  *(Mycolor4)10 &*(Mycolor5) 9 \\
      *(Mycolor1)  6 & *(Mycolor2)10 &*(Mycolor3) 9 & *(Mycolor4) 8\\
        *(Mycolor1)  8 &*(Mycolor2) 7 
\end{ytableau}$$
Then $w({\sigma}) =\color{Mycolor7}1\color{Mycolor3}2\color{Mycolor6}1\color{Mycolor5}1\color{Mycolor4}1\color{Mycolor2}2\color{Mycolor3}1\color{Mycolor1}{3}\color{Mycolor2}{14}\color{Mycolor1}2\color{Mycolor1}{1}\color{Mycolor1}4\color{Mycolor4}3\color{Mycolor3}3\color{Mycolor5}2\color{Mycolor2}3\color{Mycolor4}2$,  and
$\ch(w({\sigma})) =\maj(\sigma)  = 7. $
\end{example}

\subsection{quinv maximisers}
The following is the analogue of Proposition~\ref{lemma 1} for the $\quinv$ statistic.
\begin{proposition} \label{lemma 2} Let $\lambda=(\lambda_1 \geq \lambda_2 \geq \cdots)$ be a partition. For each $i \geq 1$, fix a set $S_{i} \subset \mathbb{N}$ of size $\lambda_{i}$. There exists a unique filling $\tau \in \mathcal{F}(\lambda)$ with  \\(a) $S_{i}  = R_i(\tau)$ for all $i$, (b) $\quinv(\tau) = n(\lambda^{\prime})$, and (c) $\maj(\tau) = \ch(w'({\tau}))$.
\end{proposition}
  \noindent
 {\em Proof:} 
Consider a filling $\tau$ with $\quinv(\tau) = n(\lambda ^{\prime})$. As before, from the definition of $\quinv$ in \eqref{eq:quinv-trip} and the maximality, we can derive the following observations about $\tau$:
\begin{enumerate}
\item For each $i$, the entries in the $i^{th}$ row of $\tau$ must be distinct.
\item For a cell $u$, if there is no cell directly below it, let $S$ be the  set consisting of $u$ and the cells to its right in $\dg(\lambda)$. Then  $\tau(u) = \min \tau(S)$.
  In particular,   $\tau$ is strictly increasing from left to right in the bottom row.  
    \item  For a cell $u$, if $v$ is the cell directly below $u$, and $S$ be the set
consisting of $u$ and the cells to its right in the same row, then: 
\begin{enumerate}
    \item if $ \tau(u) \geq \tau(v)$, then $\tau(u)$ equals $ \min \{ x\in \tau(S) : x \geq \tau(v)\},$
\item while $\tau(u) < \tau(v)$ implies that
$ x < \tau(v)$ for all $x \in\tau(S)$, and $\tau(u) = \min \tau(S)$.
\end{enumerate}
\end{enumerate}
The above facts about $\quinv$-maximum fillings implies that given a collection of sets $S_i$, if  a filling $\tau$ with the desired properties exists, then it must be unique. For the existence of $\tau$, we proceed as follows. Let $n$ denote the number of nonzero parts of $\lambda$. 
\begin{enumerate}
    \item The entries of $\tau$ in the $n^{\text{th}}$ row must be the elements of $S_n$ in strictly increasing order.
    \item Once the rows $n$ through $i + 1$ have been constructed, the entries $\tau(u)$ in the $i^{\text{th}}$ row are determined one by one, moving from left to right:
\begin{enumerate}
    \item  If there is a cell $v$ directly below $u$, then $\tau(u)$ is the smallest unused element $x$ in $S_{i}$ such that $x \geq \tau(v)$. If there is no such $x$ exists, then  $\tau(u)$ is the smallest unused element in $ S_{i}$.
    \item If there is no cell directly below $u$, then $\tau(u)$ is the smallest unused element in $S_{i}$.
\end{enumerate}
\end{enumerate}
This satisfies the properties of the first part of the proposition.

It remains to prove that $\maj(\tau) = \ch(w'({\tau}))$. Let $v_{1}$ be the first standard subword of $w({\tau})$ in the Killpatrick description of charge (section \ref{sec:killp}). Denote the occurrence of the letter $i$ in $v_{1}$ by $s_{i,{1}}$. From the definition, $s_{n,{1}}$ is the leftmost $n$ in the charge word $w'({\tau})$. From the construction (section \ref{sec:order}) of $w({\tau})$, $s_{n,{1}}$ corresponds to the smallest entry in the $n^{th}$ row of $\tau$, i.e., $\tau(n, 1)$. For $i < n$ , (i) $s_{i,{1}}$ is the leftmost $i$ to the right of $s_{i+1,1}$  in $w({\tau})$ if $\tau(i, 1) \geq \tau(i+1, 1)$ and (ii) $s_{i, {1}}$ is the leftmost $i$ otherwise.

Observe that $s_{i+1,{1}}$ lies to the right of $s_{i,{1}}$ if and only if the pair  $\tau(i+1, 1)> \tau(i, 1)$ is a descent. Therefore 
\[\ch(v_{1}) = \displaystyle\sum_{u \in \Des(C_{1}(\tau))}(\leg(u)+1)\, ,\]
where $\Des(C_{1}(\tau))$ denotes the collection of cells $u$ where the pair $u, v$ is a descent in the first column of $\tau$.

As before, if we remove the first column of $\tau$, then the remaining is a filling $\tau^{\prime}$ of $\dg(\mu)$, where $\mu$ is the shape of $\tau'$. It follows from definitions that $\tau^{\prime}$ is a $\quinv$-maximal filling of  $\dg(\mu)$. The charge word $w'({\tau^{\prime}})$ corresponding to $\tau^{\prime}$ is same as the subword obtained by removing $v_{1}$ from  $w'({\tau})$. 
Thus 
\[\ch(w'({\tau})) = \ch(v_{1}) + \ch(w'({\tau'})) \, ,\]
and 
\[ \maj(\tau) = \maj(\tau') + \displaystyle\sum_{u \in \Des(C_{1}(\tau))}(\leg(u)+1) \, . \]
We have by induction that

\[ \ch(w'({\tau})) =  \maj(\tau)\, .\] \qed

\begin{example}
\label{eg_quinm}
For the same collection of sets from Example \ref{eg_invm}, the corresponding $\quinv$-maximal filling is
$$\tau= \begin{ytableau}
      *(Mycolor1) 1 &*(Mycolor2) 2 & *(Mycolor3)7 & *(Mycolor4)3 &*(Mycolor5) 5 & *(Mycolor6) 4 & *(Mycolor7)6 \\
     *(Mycolor1)  9 &*(Mycolor2) 10 &*(Mycolor3) 7 &*(Mycolor4) 2 &*(Mycolor5) 5 \\
      *(Mycolor1)  8 &*(Mycolor2) 9 & *(Mycolor3) 6 &*(Mycolor4) 10\\
       *(Mycolor1)  7 & *(Mycolor2)8 
\end{ytableau}$$
Then $w'({\tau}) = \color{Mycolor1}1\color{Mycolor4}2\color{Mycolor2}1\color{Mycolor4}1\color{Mycolor6}1\color{Mycolor5}2\color{Mycolor5}1\color{Mycolor3}3\color{Mycolor7}1\color{Mycolor1}4\color{Mycolor3}2\color{Mycolor3}1\color{Mycolor2}4\color{Mycolor1}3\color{Mycolor2}3\color{Mycolor1}2\color{Mycolor4}3\color{Mycolor2}2$, $\ch(w'({\tau})) = 7$, and $\maj(\tau) = 7$.
\end{example}

\subsection{The first main theorem}
Given a partition $\lambda$, let $\imf(\lambda) \subset \mathcal{F}(\lambda)$ denote the set of all $\inv$-maximal fillings of $\dg(\lambda)$. Likewise, let $\qmf(\lambda)$ denote the set of  all $\quinv$-maximal fillings. Propositions \ref{lemma 1} and \ref{lemma 2}, allow us to establish the desired bijection between $\imf(\lambda)$ and $\qmf(\lambda)$. 
\begin{thm}
  There exists a bijection  $\Phi: \imf(\lambda) \to \qmf(\lambda)$ such that for all $\sigma \in \imf(\lambda)$:
   \begin{enumerate}
       \item $\sigma  \sim  \Phi(\sigma)$, 
       \item $\maj( \sigma ) = \maj( \Phi( \sigma ))$.
   \end{enumerate}
 \end{thm} 
\noindent
{\em Proof:} 
     For  $\sigma \in \imf(\lambda)$, let $R_i(\sigma)$ be the set of entries in the $i^{th}$ row of $\sigma$. From Proposition \ref {lemma 2}, there exists a unique $\tau \in  \qmf(\lambda)$ that satisfies $R_i(\tau) = R_i(\sigma)$ for all $i$. Define $\Phi(\sigma) :=  \tau$; the uniqueness of $\tau$ ensures that $\Phi$ is well defined. The combination of the existence and uniqueness assertions of Proposition \ref{lemma 1} establishes that $\Phi$ is a bijection.
     
      It is evident from the definition that $\Phi$ satisfies (1). Since $R_i(\tau) = R_i(\sigma)$ for all $i$, the ordering (as in Section \ref{sec:order})  of the cells of $\dg(\lambda)$  with respect to $\sigma$ and $<$ coincides with the ordering of the cells  with respect to $\tau$ and $<'$. Therefore, the charge words $w(\sigma)$ and $w'({\tau})$ are also the same, and Propositions \ref{lemma 1} and  \ref{lemma 2}, together imply that $\maj(\sigma) = \ch(w({\sigma})) = \ch(w'({\tau})) = \maj(\tau)$. \qed

\section{The bijection in the modifed Hall-Littlewood case}
Theorem \ref{Theorem 1}, Theorem \ref{Theorem 2}, and \eqref{eq1.1} together imply that the modified Hall-Littlewood polynomial can be expressed as
\begin{equation}
    \widetilde{H}_{\lambda}(x;t) = \sum_{\substack{\sigma \in \mathcal{F}(\lambda) \\ \inv(\sigma)=0}}t^{\maj(\sigma)} x^{\sigma} \, ,
\end{equation}
or 
\begin{equation}
    \widetilde{H}_{\lambda}(x;t) = \sum_{\substack{\sigma \in \mathcal{F}(\lambda)\\ \quinv(\sigma)=0}}t^{\maj(\sigma)} x^{\sigma} \, .
\end{equation}
The following two propositions characterise the class of $\inv$ zero and $\quinv$ zero fillings of $\dg(\lambda)$.
\subsection{inv-zero fillings}
\begin{proposition}
\label{L.3} Let $\lambda=(\lambda_1 \geq \lambda_2 \geq \cdots)$ be a partition.
 For each $i \geq 1$, let $M_{i}$ be a multiset of positive integers, with size $|M_{i}| = \lambda_{i}$. Then there exists a unique filling $\sigma \in \mathcal{F}(\lambda)$ such that  $\inv(\sigma) = 0$ and $R_i(\sigma) = M_{i}$ for all $i$. Further, $ \sigma $ satisfies $\maj( \sigma ) = \cch(cw({\sigma}))$.
\end{proposition}
The proof of Proposition \ref{L.3} can be seen in \cite[Proposition 7.1]{HHL01}.
\begin{example}
\label{ex_invz}
Let $\lambda = (7, 5, 4, 2) $, and let $M_{1}= \{1, 2, 2, 3, 4, 4, 4\}, \, M_{2}= \{1, 2, 3, 4, 4\},\\\,M_{3}= \{2, 2, 4, 5\},\,\text{and } M_{4}= \{3, 5\}$, then the corresponding $\inv$ zero filling is $$\sigma = \begin{ytableau}
      *(Mycolor1) 1 &*(Mycolor2) 2 &*(Mycolor3) 2 & *(Mycolor4)3 & *(Mycolor5)4 &*(Mycolor6)4 &*(Mycolor7) 4 \\
      *(Mycolor1) 2 & *(Mycolor2)3 &*(Mycolor3) 4 & *(Mycolor4)4 & *(Mycolor5)1 \\
      *(Mycolor1)  4 & *(Mycolor2)5 &*(Mycolor3)2 & *(Mycolor4)2\\
       *(Mycolor1)  5 & *(Mycolor2)3 
\end{ytableau}$$
Then $cw({\sigma}) =  \color{Mycolor2}3\color{Mycolor1}4\color{Mycolor7}1\color{Mycolor6}1\color{Mycolor5}1\color{Mycolor4}2\color{Mycolor3}2\color{Mycolor1}3\color{Mycolor4}1\color{Mycolor2}2\color{Mycolor2}4\color{Mycolor3}1\color{Mycolor2}1\color{Mycolor1}2\color{Mycolor4}3\color{Mycolor3}3\color{Mycolor1}1\color{cyan}2$,  and
$ \maj(\sigma) = \cch( cw({\sigma})) = 15 .$
\end{example}
\subsection{quinv-zero fillings}
\noindent
The proposition below is the analogue of Proposition \ref{L.3} for $\quinv$ zero fillings.
\begin{proposition}
\label{L.4} Let $\lambda=(\lambda_1 \geq \lambda_2 \geq \cdots)$ be a partition.
 For each $i \geq 1$, let $M_{i}$ be a multiset of positive integers with size $|M_{i}| = \lambda_{i}$. Then there exists a unique filling $\tau \in \mathcal{F}(\lambda)$ such that $\quinv(\tau) = 0$ and $R_i(\tau) = M_{i}$ for all $i$. Further, $\tau$ satisfies 
  $\maj(\tau) = \cch(cw'({\tau})).$
\end{proposition}

\noindent
  {\em Proof:} 
 We observe the following facts about a $\quinv$ zero filling $\tau \in \mathcal{F}(\lambda)$:
 \begin{enumerate}
     \item Let $u$ be a cell in $\dg(\lambda)$  which does not have a cell below, and let $S$ be the set consisting of $u$ and the cells to its right in the same row. Then $\tau(u) = \max\{\tau(w) : w \in S\}$, in particular $\tau$ is non-increasing in the bottom row. 
     \item For a cell $u$, if there  exists  a cell $v$ directly below it, and let $S'$ be the set consisting of $u$ and the cells to its right in the same row, then
\begin{enumerate}
    \item  if $\tau(u) \geq \tau(v)$, then $\tau(v) \leq x$ for all $x \in \tau(S')$ and $\tau(u)= \text{ max } \{ \tau(w) : w \in S' \}$, or
    \item if $\tau(u) < \tau(v) $ , then $\tau(u) = \text{ max } \{ \tau(w) \text{ for all }  w \in S' : \tau(w) < \tau(v)  \} $.
\end{enumerate}
\end{enumerate}
From the above observations, it is clear that if there exists a  $\quinv$ zero filling corresponding to a given collection $\text{M} = \{M_{i} : M_{i} \text{ is a multiset }, |M_{i}| = \lambda_{i}, i \geq 1\}$, then it must be unique.

   For the existence, we construct a filling $\tau $, with  $\quinv(\tau) =0$, and  the multiset of entries in the $i^{th}$ row of $\tau$, $M_{\tau_{i}} = M_{i}$ for each $i$.
  Let $n$ denote the number of nonzero parts of $\lambda$.  Since  $\quinv(\tau) =0 $, $\tau$ must be non increasing in the $n^{th}$ row, therefore its entries must be the elements of $M_{n}$ in non-increasing order. After constructing rows $n$ through $i+1$, the entries $\tau(u)$ in row $i$ are determined sequentially from left to right as follows 
   \begin{enumerate}
       \item If there is a cell $v$  directly below $u$, and if $M_{i}$ contains an unused element $ x < \tau(v)$, then $\tau(u)$ is defined to be the largest such $x$; otherwise $\tau(u)$ is assigned the largest unused $x \in M_{i}$.
       \item If there is no cell directly below $u$, then $\tau(u)$ is assigned the largest unused $x \in M_{i}$.
   \end{enumerate}
This completes the the first part of the proof. 

To prove $\maj(\tau) = \cch(cw'({\tau}))$, let $cw_{1}'$ be the first standard subword of $cw'({\tau})$ in the Killpatrick description of charge (or cocharge). Denote each alphabet $i$ in $cw_{1}'$ as $r_{i,{1}}$. From the definition, $r_{n,{1}}$ is the leftmost $n$ in $cw'({\tau})$, and from the construction of $cw'({\tau})$, $r_{n,{1}}$ corresponds to the largest entry in the last row of $\tau$, that is $r_{n,{1}}$ corresponds $\tau(n, 1)$.
For $i < n$ , $r_{i,{1}}$ will be (i) the leftmost $i$ to the right of $r_{{i+1},{1}}$ in $cw'({\tau})$ if $\tau(i, 1) \geq \tau(i+1, 1)$, and (ii) the rightmost $i$ otherwise. And observe that $r_{i+1, {1}}$ lies in the right of $r_{i,{1}}$ in $cw_{1}'$ if and only if the pair  $\tau(i+1, 1)> \tau(i, 1)$ is a descent.
Therefore 
\[\cch(cw_{1}') = \displaystyle\sum_{u \in \text{Des}(C_{1}(\tau))}(\text{leg}(u)+1)\, ,\] 
where $\text{Des}(C_{1}(\tau))$ denotes the collection of cells $u$ where the pair $\tau(u)> \tau(v)$ is a descent in the first column of $\tau$.

If we remove the first column of $\tau$, then the remaining will be a $\quinv$-zero filling $\tau^{\prime}$ of $\dg(\mu)$, where $\mu$ is the shape of $\tau'$. The {cocharge word } $cw'({\tau'})$ corresponding to $\tau'$ is the same as the word obtained by removing $cw_{1}'$ from  $cw'({\tau})$. Thus 
\[\cch(cw'({\tau})) = \cch(cw_{1}') + \cch(cw'({\tau'}))\]
and 
\[ \maj(\tau) =  \displaystyle\sum_{u \in \Des(C_{1}(\tau))}(\leg(u)+1) + \maj(\tau') \, .\]
By induction we have
\[ \cch(cw'({\tau})) =  \maj(\tau) \, .\] \qed

\noindent
\begin{example}
For the same collection of multisets from Example \ref{ex_invz}, the unique $\quinv$ zero filling is
$$\tau = \begin{ytableau}
     *(Mycolor1)  2 & *(Mycolor2)4 &*(Mycolor3) 3 &*(Mycolor4) 2 &*(Mycolor5) 1 &*(Mycolor6) 4 &*(Mycolor7) 4 \\
      *(Mycolor1) 3 &*(Mycolor2) 1 &*(Mycolor3) 4 &*(Mycolor4) 4 &*(Mycolor5) 2 \\
       *(Mycolor1) 4 & *(Mycolor2)2 &*(Mycolor3) 5 & *(Mycolor4)2\\
       *(Mycolor1)  5 &*(Mycolor2) 3 
\end{ytableau}$$
 Then  $cw'({\tau}) = \color{Mycolor3} 3\color{Mycolor1}4\color{Mycolor2}1\color{Mycolor6}1\color{Mycolor7}1\color{Mycolor3}2\color{Mycolor4}2\color{Mycolor1}3\color{Mycolor3}1\color{Mycolor1}2\color{Mycolor2}4\color{Mycolor1}1\color{red}1\color{Mycolor5}2\color{Mycolor2}3\color{Mycolor4}3\color{cyan}1\color{Mycolor2}2$, and $ \maj(\tau) = \cch( cw'({\tau}))= 15$.
\end{example}

\subsection{The second main theorem }
 Given a partition $\lambda$, let $\izf(\lambda)$ (resp. $\qzf(\lambda)$) denote the set of $\inv$-zero (resp. $\quinv$-zero) fillings of $\dg(\lambda)$. Propositions \ref{L.3} and \ref{L.4}, now give us the desired bijection between $\izf(\lambda)$ and $\qzf(\lambda)$. 
\noindent
\begin{thm}
  Let $\lambda$ be a partition. There exist a bijection  $\varphi$  from $\izf(\lambda)$  to $\qzf(\lambda)$ such that for all $\sigma \in \izf(\lambda)$:
\begin{enumerate}
\item $\sigma  \sim  \varphi(\sigma)\, ,$ 
\item $\maj( \sigma ) = \maj( \varphi( \sigma ))\, .$
\end{enumerate}
\end{thm}

{\em Proof:}
 Consider $\sigma \in \izf(\lambda)$, and let $M_{\sigma_{i}}$ be the multiset of entries in the $i^{th}$ row of $\sigma$. From Proposition \ref{L.4} there exist a unique filling $\tau $ in $ \qzf(\lambda)$, such that $M_{\sigma_{i}} = M_{\tau_{i}}$. Define $\varphi(\sigma) := \tau $, then $\varphi$ is well defined, and Proposition \ref{L.3} shows that it is in fact a bijection.
 
 It is evident from the definition that $\varphi$ satisfies (1).
 Since $M_{\sigma_{i}} = M_{\tau_{i}}$ for all $i$, the orderings of the cells of $\dg(\lambda)$ with respect to $\sigma$ and $<$ and  the ordering with respect to $\tau$ and $<'$ must be same, so the cocharge words $cw({\sigma})$ and  $cw'({\tau})$ must coincide. Therefore, the Propositions \ref{L.3} and \ref{L.4} together imply that
 $\maj(\sigma) = \cch(cw({\sigma})) = \cch(cw'({\tau})) =  \maj(\tau)$. \qed

\bibliography{qw_mhl}
\end{document}